\documentclass[12pt]{amsart}
\usepackage{amsfonts}
\usepackage{amssymb}
\numberwithin{equation}{section}
\theoremstyle{plain}
 \newtheorem{thm}{Theorem}[section]
 \newtheorem{lem}{Lemma}[section]
 \newtheorem{cor}{Corollary}[section] 
 \newtheorem{prop}{Proposition}[section]
\theoremstyle{definition}
 \newtheorem{defn}{Definition}[section]

\newcommand{\mcal}{\mathcal}

\newcommand{\q}{\quad}

\newcommand{\R}{\mathbb{R}}
\newcommand{\N}{\mathbb{N}}

\begin{document}
\setlength{\baselineskip}{15pt}
\setlength{\parindent}{1.8pc}

\title{Two  non-closure properties on the class of subexponential densities}
\author{Toshiro Watanabe and Kouji Yamamuro}
\maketitle
\footnote[0]{T. Watanabe: Center for Mathematical Sciences The Univ. of Aizu. Aizu-Wakamatsu 965-8580, Japan. e-mail: t-watanb@u-aizu.ac.jp\\ 
 K. Yamamuro: Faculty of Engineering Gifu Univ. Gifu 501-1193, Japan. \\ e-mail: yamamuro@gifu-u.ac.jp}\\

{\small
{\bf Abstract}.   Relations between  subexponential densities and  locally subexponential distributions are discussed. It is shown that the class of subexponential densities is neither closed under convolution roots nor closed  under asymptotic equivalence. A remark is  given on the closure under convolution roots  for the class of convolution equivalent  distributions.}\\

\medskip  

{\bf Key words or phrases: } subexponential densities, local subexponentiality, convolution roots, asymptotic equivalence\\
             
{\bf Mathematics Subject Classification}: 60E99, 60G50   \\

\section{Introduction and main results}
In what follows, we denote by  $\R$ the real line and by  $\R_+$ the half line $[0,\infty)$. Let $\N$ be the totality of positive integers. The symbol $\delta_a(dx)$ stands for the delta measure at $a\in \R$. Let $\eta$ and $\rho$ be probability measures on $\R$. We denote the convolution of $\eta$ and $\rho$ by $\eta*\rho$ and denote $n$-th convolution power of $\rho$ by  $\rho^{n*}$.  
 Let $f(x)$ and $g(x)$ be integrable functions on $\R$. We denote by $f^{n\otimes}(x)$ $n$-th convolution power of  $f(x)$ and by $f\otimes g(x)$ the convolution of $f(x)$ and $g(x)$.  For positive  functions $f_1(x)$ and $g_1(x)$ on $[a,\infty)$ for some $a \in \R$, we define the relation $f_1(x) \sim g_1(x)$ by $\lim_{x \to \infty}f_1(x)/g_1(x) =1$. We also define the relation $a_n \sim b_n$ for positive sequences $\{a_n\}_{n=A}^{\infty}$ and $\{b_n\}_{n=A}^{\infty}$ with $A \in \N$ by $\lim_{n \to \infty}a_n/b_n =1$. We define the class  ${\mcal P}_+$ as the totality of probability distributions on $\R_+$. In this paper, we prove that the class of subexponential densities is not closed under two important closure properties. We say that a measurable function  $g(x)$  on $\R$ is a density function if $\int_{-\infty}^\infty g(x)dx=1$ and $g(x)\geq 0$ for all $x\in\R$.

\begin{defn}
(i) A nonnegative measurable function  $g(x)$  on $\R$  belongs to the class ${\bf L}$ if $g(x)>0$ for all sufficiently large $x>0$ and if $g(x+a) \sim g(x)$ for any $a \in \R$.

(ii)  A measurable function  $g(x)$  on $\R$  belongs to the class ${\mcal L}_{d}$  if $g(x)$ is a density function and $g(x)\in{\bf L}$.

(iii) A measurable function  $g(x)$  on $\R$  belongs to the class ${\mcal S}_{d}$  if $g(x)\in{\mcal L}_d$ and $g\otimes g(x) \sim 2g(x)$.

(iv)  A distribution $\rho$  on $\R$   belongs to the class ${\mcal L}_{ac}$ if there is $g(x)\in {\mcal L}_d$ such that $\rho(dx) =g(x) dx$. 

(v) A distribution $\rho$   on $\R$   belongs to the class ${\mcal S}_{ac}$  if  there is $g(x)\in {\mcal S}_d$ such that  $\rho(dx) =g(x) dx$.

\end{defn}
Densities in the class ${\mcal S}_{d}$ are called {\it subexponential densities} and those in the class ${\mcal L}_{d}$ are called {\it long-tailed densities}. The study on the class ${\mcal S}_{d}$ goes back to Chover et al.\ \cite{cnw}. Let $\rho$ be a  distribution on $\R$. Note that $c^{-1}\rho((x-c,x])$ is a density function  on $\R$ for every $c>0$.

\begin{defn}
(i)   Let $\Delta:= (0,c]$ with $c>0.$ A distribution $\rho$    on $\R$   belongs to the class $ {\mcal L}_{\Delta}$   if $\rho((x,x+c])  \in   {\bf L}$.
 
(ii)  Let $\Delta:= (0,c]$ with $c>0.$  A distribution $\rho$    on $\R$   belongs to the class $ {\mcal S}_{\Delta}$   if $\rho  \in  {\mcal L}_{\Delta}$ and 
$\rho*\rho((x,x+c])  \sim 2 \rho((x,x+c]). $ 

(iii) A distribution  $\rho$    on $\R$   belongs to the class ${\mcal L}_{loc}$   if $\rho  \in {\mcal L}_{\Delta}$ for each  $\Delta:= (0,c]$ with $c>0.$ 

(iv)  A distribution  $\rho$   on $\R$   belongs to the class ${\mcal S}_{loc}$  if $\rho  \in {\mcal S}_{\Delta}$ for each $\Delta:= (0,c]$ with $c>0.$

(v) A distribution $\rho \in \mathcal{L}_{loc}$ belongs to the class $\mcal{ UL}_{loc}$  if there exists $p(x) \in  \mathcal{L}_{d}$ such that 
$c^{-1}\rho((x-c,x]) \sim p(x)$ uniformly in $c \in (0,1]$.

(vi) A distribution $\rho \in \mathcal{S}_{loc}$ belongs to the class $\mcal{ US}_{loc}$  if there exists $p(x) \in  \mathcal{S}_{d}$ such that 
$c^{-1}\rho((x-c,x]) \sim p(x)$ uniformly in $c \in (0,1]$.
\end{defn}

Distributions in the class ${\mcal S}_{loc}$ are called {\it locally subexponential}, those in the class ${\mcal US}_{loc}$ are called {\it uniformly locally subexponential}. The class $ {\mcal S}_{\Delta}$ was introduced by Asmussen et al.\ \cite{afk03} and the class ${\mcal S}_{loc}$ was by Watanabe and Yamamuro \cite{wy10}.  Detailed acounts of the classes ${\mcal S}_{d}$ and ${\mcal S}_{\Delta}$ are found in the book of Foss et al.\ \cite{fkz13}.
First, we present some interesting results on the classes ${\mcal S}_{d}$ and ${\mcal S}_{loc}$. 

\begin{prop} We have the following.

(i) Let $\Delta:= (0,c]$ with $c>0$ and let $p(x):=c^{-1}\mu((x-c,x])$ for a distribution $\mu$ on  $\R_+$.   Then $\mu \in  {\mcal S}_{\Delta}$ if and only if $p(x)\in {\mcal S}_{d}$. Moreover,  $\mu \in {\mcal S}_{loc}\cap {\mcal P}_+$    if and only if there exists a density function $q(x)$ on $\R_+ $ such that $q(x)\in\mathcal{S}_{d}$ and 
$c^{-1}\mu((x-c,x])\sim q(x)$ for every $c >0$.

(ii) Let $\rho_1(dx):=q_1(x)dx$ be a distribution on $\R_+$. If $q_1(x)$ is continuous with compact support and if $\rho_2 \in {\mcal S}_{loc}\cap {\mcal P}_+$, then $\rho_1*\rho_2(dx)=\left(\int_{0-}^{x+}q_1(x-u)\rho_2(du)\right)dx$ and $ \int_{0-}^{x+}q_1(x-u)\rho_2(du)\in {\mcal S}_{d}$.

(iii) Let $\mu$ be a distribution on $\R_+$. If there exist   distributions $\rho_c$ for $c>0$ such that, for every $c>0$,  the support of $\rho_c$ is included in $[0,c]$ and $\rho_c*\mu \in \mathcal{S}_{loc}$, then $\mu\in \mathcal{S}_{loc}$. 
\end{prop}

\begin{defn}
(i) We say that a class $\mathcal{C}$ of probability distributions  on $\R$ is closed under convolution roots if  $\mu^{n*} \in \mathcal{C}$ for some $ n \in \mathbb{N}$ implies that $\mu \in \mathcal{C}$. 

(ii) Let  $p_1(x)$ and  $p_2(x)$ be density functions on $\R$.   We say that a class $\mathcal{C}$ of density functions is closed under asymptotic equivalence if  $p_1(x) \in \mathcal{C}$ and
$p_2(x) \sim c p_1(x)$ with $c >0$ implies that $p_2(x)\in \mathcal{C}$.
\end{defn}

The class $\mathcal{S}_{ac}$ is a proper subclass of the class $\mathcal{US}_{loc}$ because a distribution in $\mathcal{US}_{loc}$ can have a point mass. Moreover,  the class $\mathcal{US}_{loc}$ is a proper subclass of the class $ \mathcal{S}_{loc}$ as the following theorem shows.
 
\begin{thm}

 There  exists a distribution $\mu \in \mathcal{S}_{loc}\setminus \mathcal{US}_{loc}$ such that $\mu^{2*} \in {\mcal S}_{ac}$.

\end{thm} 

\begin{cor} We have the following.

(i) The class $\mathcal{S}_{ac}$ is not closed under convolution roots.

(ii) The class $\mathcal{US}_{loc}$ is not closed under convolution roots.

(iii)  The class $\mathcal{L}_{ac}$ is not closed under convolution roots.

(iv) The class $\mathcal{UL}_{loc}$ is not closed under convolution roots.
\end{cor}

The class $\mathcal{S}_{d}$ is closed under asymptotic equivalence in the one-sided case. See (ii) of Lemma 2.1 below. However, Foss et al.\ \cite{fkz13} suggest the possibility of non-closure under asymptotic equivalence for the class $\mathcal{S}_{d}$ in the two-sided case. We exactly prove it as follows.

\begin{thm} 
  The class $\mathcal{S}_{d}$ is not closed under asymptotic equivalence, that is, there exist $p_1(x)\in \mathcal{S}_{d}$ and  $p_2(x) \notin \mathcal{S}_{d}$  such that $p_2(x) \sim c p_1(x)$ with $c >0$. 
\end{thm}

In Sect.\ 2, we prove Proposition 1.1. In Sect.\ 3, we prove Theorems 1.1 and 1.2.  In Sect.\ 4, we give a remark on the closure under convolution roots. 


\section{Proof of Proposition 1.1}
We present two lemmas for the proofs of  main results and then prove Proposition 1.1.

\begin{lem} Let $f(x)$ and $g(x)$ be density functions  on $\R_+$.

(i) If $f(x)\in {\mcal L}_{d}$, then $f^{n\otimes}(x) \in {\mcal L}_{d}$ for every $n \in \N$.

(ii) If $f(x) \in {\mcal S}_{d}$ and $g(x) \sim c f(x)$ with $c >0$, then  $g(x) \in {\mcal S}_{d}$.

(iii) Assume that $f(x) \in {\mcal L}_{d}$. Then, $f(x)\in {\mcal S}_{d}$ if and only if
\begin{equation} 
\lim_{A \to \infty} \limsup_{x \to \infty} \frac{1}{f(x)}\int_{A}^{x-A}f(x-u)f(u)du =0.\nonumber
\end{equation}
\end{lem}

{\it Proof }  \q  Proof of assertion (i) is due to Theorem 4.3 of \cite{fkz13}. Proofs of assertions (ii) and (iii) are due to Theorems 4.8 and 4.7 of \cite{fkz13}, respectively.\qed

\begin{lem}
(i) Let $\Delta:= (0,c]$ with $c>0.$ Assume that $\rho \in  {\mcal L}_{\Delta} \cap{\mcal P}_+$. Then, $\rho \in  {\mcal S}_{\Delta}$ if and only if 
\begin{equation} 
\lim_{A \to \infty} \limsup_{x \to \infty} \frac{1}{\rho((x,x+c])}\int_{A+}^{(x-A)-}\rho((x-u,x+c-u])\rho(du) =0.\nonumber
\end{equation}

(ii) Assume that $\rho \in {\mcal L}_{loc}\cap{\mcal P}_+$.  
Then, $\rho^{n*} \in {\mcal L}_{loc}$ for every $n \in \N$. Moreover, $\rho((x-c,x]) \sim c\rho((x-1,x])$ for every $c>0$. 

(iii) Let $\rho_2\in{\mcal P}_+$. If $\rho_1 \in {\mcal S}_{loc}\cap {\mcal P}_+$ and $\rho_2((x-c,x]) \sim c_1 \rho_1((x-c,x]) $ with $c_1 >0$ for every $c >0$, then  $\rho_2 \in {\mcal S}_{loc}\cap {\mcal P}_+$.

\end{lem}
{\it Proof } \q   Proof of assertion (i) is due to Theorem 4.21 of  \cite{fkz13}. First assertion of (ii) is due to Corollary 4.19 of \cite{fkz13}. Second one  is proved as (2.6) in Theorem 2.1 of \cite{wy10}. Proof of assertion (iii) is due to Theorem 4.22 of \cite{fkz13}.\qed

\medskip

{\it Proof of (i) of Proposition 1.1 } \q 
Let $\rho(dx):=c^{-1}1_{[0,c)}(x)dx$. First, we prove that if $\mu \in {\mcal S}_{loc}\cap{\mcal P}_+$, then $\rho*\mu \in {\mcal S}_{ac}$. We can assume that $c=1$. Suppose that $\mu \in {\mcal S}_{loc}$. Let $p(x):=  \mu((x-1,x])$.  We have $\rho*\mu(dx)=   \mu((x-1,x])dx$ and hence $p(x) \in {\mcal L}_{d}$. Let $A$ be a positive integer and let $X, Y$ be independent random variables with the same distribution $\mu$.  Then, we have for $x > 2A+2$ 
\begin{equation}
\begin{split}
  & \int_A^{x-A}p(x-u)p(u)du\\ \nonumber
 &  =2\int_A^{x/2}p(x-u)p(u)du\\ \nonumber
  & = 2\int_A^{x/2}P( x-u-1<X \le x-u, u-1<Y \le u)du\\ \nonumber
 & \le 2\int_A^{x/2}P(X>A, Y >A, x-2 <X+Y \le x, u-1<Y \le u)du\\ \nonumber
 & \le 2\sum_{n=A}^{\infty}\int_n^{n+1}P(X>A, Y >A, x-2 <X+Y \le x, n-1<Y \le n+1)du\\ \nonumber
 & \le 4P(X>A, Y >A, x-2 <X+Y \le x)\\ \nonumber
& \le 4\int_{A+}^{(x-A)-} \mu((x-2-u,x-u])\mu(du). \nonumber
\end{split}
\end{equation}
Since  $\mu \in {\mcal S}_{loc}$, we obtain from  (i) of Lemma 2.2 that
$$ \lim_{A \to \infty}\limsup_{x \to \infty}\frac{\int_A^{x-A}p(x-u)p(u)du }{p(x)}=0.$$
Thus, we see from (iii) of Lemma 2.1 that $p(x)\in {\mcal S}_{d}$. 

Conversely, suppose that $p(x)\in {\mcal S}_{d}$. Then, we have $\mu \in {\mcal L}_{\Delta}$. Let $[y]$ be the largest integer not exceeding a real number $y$.  Choose sufficiently large integer $A > 0$. Note that there are positive constants $c_j$ for $1 \le j \le 4$ such that  
\begin{equation}
c_1p(x-n)\le p(x-u) \le  c_2p(x-n)\mbox{ and }c_3p(n)\le p(u) \le  c_4p(n) \nonumber
\end{equation}
 for $n \le u \le n+1$, $A \le n \le [x+1-A]$, and $x > 2A+2$. Thus, we find that
\begin{equation}
\begin{split}
& P(A < X, A < Y, x < X+Y \le x+1)\\ \nonumber
& \le \sum_{n=A}^{[x+1-A]}\int_n^{n+1}\mu((x-u,x+1-u])\mu(du)\\ \nonumber
& =\sum_{n=A}^{[x+1-A]}\int_n^{n+1}p(x-u+1)\mu(du)\\ \nonumber
& \le c_2\sum_{n=A}^{[x+1-A]}p(x-n+1)p(n+1)\\ \nonumber
& \le \frac{c_2}{c_1c_3}\sum_{n=A}^{[x+1-A]}\int_n^{n+1}p(x-u+1)p(u+1)du\\ \nonumber
& \le \frac{c_2}{c_1c_3}\int_A^{x+2-A}p(x-u+1)p(u+1)du\\ \nonumber
\end{split}
\end{equation}
Since $p(x)\in {\mcal S}_{d}$, we establish from (iii) of Lemma 2.1 that 
$$ \lim_{A \to \infty}\limsup_{x \to \infty}\frac{P(A < X, A < Y, x < X+Y \le x+1) }{P(x< X \le  x+1)}=0.$$
Thus,  $\mu \in  {\mcal S}_{\Delta}$ by (i) of Lemma 2.2. Note from (ii) of Lemma 2.2 that if $\mu \in  {\mcal S}_{loc}$, then $c^{-1}\mu((x-c,x])\sim \mu((x-1,x])$ for every $c >0$. Thus, the second assertion is true. \qed

\medskip 

{\it Proof of (ii) of Proposition 1.1 } \q 
Suppose that $\rho_1(dx):=q_1(x)dx$ be a distribution on $\R_+$ such that $q_1(x)$ is continuous with compact support in $[0,N]$.
 Let $q(x):=\int_{0-}^{x+}q_1(x-u)\rho_2(du)$. 
For  $ M \in \N$, there are $\delta(M) >0$ and $a_n =a_n(M) \ge 0$ for $n \in \N$ such that $\lim_{M \to \infty}\delta(M) =0$ and $a_n \le q_1(x) \le a_n +\delta(M)$ for $M^{-1}(n-1) < x \le M^{-1}n$  and $ 1 \le n  \le MN$.   Define $J(M; x)$ as
\begin{equation}
J(M; x):= \sum_{n=1}^{MN}a_n \rho_2((x-M^{-1}n, x-M^{-1}(n-1)]).  \nonumber
\end{equation}
Then, we have 
\begin{equation}
J(M; x) \sim \rho_2((x-1, x])\sum_{n=1}^{MN}a_n M^{-1} 
\end{equation}
and for $x > N$
\begin{equation}
J(M; x) \le q(x) \le J(M; x) +\delta(M) \rho_2((x-N, x]). \nonumber
\end{equation}
Since $\lim_{M \to \infty}\delta(M) =0$ and
\begin{equation}
\lim_{M \to \infty}\sum_{n=1}^{MN}a_n M^{-1} = \int_0^N q_1(x)dx =1, \nonumber
\end{equation}
we obtain from (2.1) that 
\begin{equation}
q(x) \sim \rho_2((x-1, x]). \nonumber
\end{equation}
Since $\rho_2 \in {\mcal S}_{loc}$, we conclude from (i) of Proposition 1.1 that
$q(x)\in {\mcal S}_{d}.$  \qed

\medskip 

{\it Proof of (iii) of Proposition 1.1 } \q  Suppose that the support of $\rho_c$ is included in $[0,c]$ and $\rho_c*\mu \in \mathcal{S}_{loc}$ for every $c>0$. Let $X$ and $Y$ be independent random variables with the same distribution $\mu$, and let $X_c$ and $Y_c$ be independent random variables with the same distribution $\rho_c$. 
Define $J_1(c;c_1;a;x)$ and $J_2(c;c_1;a;x)$ for $a \in \R$ and $c_1> 0$ as
\begin{equation}
J_1(c;c_1;a;x):=\frac{P( x +a< X +X_c\le x+c_1+a)}{ P( x < X +X_c\le x+c_1+c)}, \nonumber
 \end{equation}
\begin{equation} 
J_2(c;c_1;a;x):=\frac{P( x +a< X +X_c\le x+c_1+c+a)}{ P( x < X +X_c\le x+c_1)}. \nonumber
\end{equation}
We see that
\begin{equation}
J_1(c;c_1;a;x)\le \frac{P( x +a< X\le x+c_1+a)}{ P( x < X\le x+c_1)}\le J_2(c;c_1;a;x).
\end{equation}
Since $\rho_c*\mu \in \mathcal{L}_{loc}$, we obtain that 
\begin{equation}
\lim_{x \to \infty} J_1(c;c_1;a;x)= \frac{c_1}{c_1 +c} \nonumber
\end{equation}
and 
\begin{equation}
   \lim_{x \to \infty} J_2(c;c_1;a;x)= \frac{c_1+c}{c_1 }. \nonumber
\end{equation}
Thus, as $c \to 0$ we have by (2.2)
\begin{equation}
\lim_{x \to \infty}\frac{P( x +a< X\le x+c_1+a)}{ P( x < X\le x+c_1)}=1, \nonumber
\end{equation}
and hence $\mu \in \mathcal{L}_{loc}$. We find from $\rho_c*\mu \in \mathcal{S}_{loc}$ and (i) of Lemma 2.2 that 
\begin{equation}
\begin{split}
&\lim_{A \to \infty}\limsup_{x \to \infty}\frac{P(X>A, Y>A,  x< X+Y\le x+c_1)}{ P( x < X\le x+c_1)}\\ \nonumber
&\le \lim_{A \to \infty}\limsup_{x \to \infty}\frac{P(X>A, Y>A,  x< X+X_c+Y+Y_c\le x+c_1+2c)}{ P( x < X+X_c\le x+c_1)}=0. \nonumber
\end{split}
\end{equation}
Thus, we see from (i) of Lemma 2.2 that $\mu \in \mathcal{S}_{loc}$. \qed


\medskip
\section{Proofs of Theorems 1.1 and 1.2}

For the proofs of the theorems, we introduce a distribution $\mu$ as follows.
Let $1 <x_0 < b$ and choose $\delta \in (0,1)$ satisfying $\delta< (x_0-1)\land(b-x_0)$. We take a continuous periodic function $h(x)$ on $\mathbb{R}$ with period $\log b$ such that $h(\log x) >0$ for $x \in [1,x_0)\cup (x_0,b]$ and
\begin{eqnarray*}
h(\log x) &=&\left\{
\begin{array}{ll}
0  & \mbox{for $x=x_0$},\\
{\displaystyle\frac{-1}{\log |x-x_0|}} & \mbox{for each $x$ with $0<|x-x_0|<\delta$}.
\end{array}
\right.
\end{eqnarray*}
Let 
\begin{equation}
\phi(x):=x^{-\alpha-1}h(\log x) 1_{[1,\infty)}(x)\nonumber
\end{equation}
 with $\alpha >0$. Here, the symbol $1_{[1,\infty)}(x)$ stands for the indicator function of the set $[1,\infty)$. Define a distribution $\mu$ as
\begin{eqnarray*}
&& \mu(dx):=M^{-1}\phi(x)dx,
\end{eqnarray*}
where $M:=\int_1^\infty x^{-1-\alpha}h(\log x)dx$.

\vspace{2ex}

\begin{lem}
We have $\mu\in\mathcal{L}_{loc}$.
\end{lem}

{\it Proof}\quad Let $\{y_n\}$ be a sequence such that $1\leq y_n\leq b$ and $\lim_{n\to\infty}y_n=y$ for some $y\in[1,b]$. Then, we put $x_n=b^{m_n}y_n$, where $m_n$ is a positive integer and $\lim_{n\to\infty}x_n=\infty$. In what follows, $c>0$ and $c_1\geq 0$. \par
Case 1. \quad Suppose that $y\not=x_0$. Let $x_n+c_1\leq u\leq x_n+c_1+c$. Then, we have
\begin{eqnarray}
&& y_n+b^{-m_n}c_1\leq b^{-m_n}u\leq y_n+b^{-m_n}(c_1+c),\label{lem0100}
\end{eqnarray} 
and thereby $\lim_{n\to\infty}b^{-m_n}u=y$. This yields that
\begin{eqnarray*}
h(\log u)=h(\log(b^{-m_n}u))\sim h(\log y).
\end{eqnarray*}
Hence, we obtain that
\begin{eqnarray*}
\int_{x_n+c_1}^{x_n+c_1+c}\phi(u)du &=& \int_{x_n+c_1}^{x_n+c_1+c}u^{-1-\alpha}h(\log u)du\\
&\sim& x_n^{-1-\alpha}\int_{x_n+c_1}^{x_n+c_1+c}h(\log u)du\sim c x_n^{-1-\alpha}h(\log y),
\end{eqnarray*}
so that
\begin{eqnarray}
\int_{x_n}^{x_n+c}\phi(u)du &\sim& \int_{x_n+c_1}^{x_n+c_1+c}\phi(u)du \label{lem0101}
\end{eqnarray}

\par

Case 2.\quad Suppose that $y=x_0$.  Let $x_n+c_1\leq u\leq x_n+c_1+c$ and put
\begin{eqnarray*}
E_n:=\{u\ :\ |b^{-m_n}u-x_0|\leq \epsilon b^{-m_n}\},
\end{eqnarray*}
where $\epsilon>0$. For sufficiently large $n$, we have for $u\in E_n$
\begin{eqnarray}
&& -\log|b^{-m_n}u-x_0|\geq -\log\epsilon b^{-m_n}\geq \frac{1}{2}m_n\log b\label{lem0102}
\end{eqnarray}
Set $\lambda_n:=|y_n-x_0|b^{m_n}$. It suffices that we consider the case where there exists a limit of $\lambda_n$ as $n\to\infty$, so we may put $\lambda:=\lim_{n\to\infty}\lambda_n$. This limit permits infinity. We divide $\lambda$ in the two cases where $\lambda<\infty$  and $\lambda=\infty$.  \par
Case 2-1. Suppose that $0\leq \lambda<\infty$. Now, we have
\begin{eqnarray*}
 && \int_{x_n+c_1}^{x_n+c_1+c}h(\log u)du\\
 &=& \int_{[x_n+c_1, x_n+c_1+c]\backslash E_n}h(\log u)du+\int_{[x_n+c_1, x_n+c_1+c]\cap E_n}h(\log u)du.
\end{eqnarray*}
Let $u\in[x_n+c_1, x_n+c_1+c]\backslash E_n$. For sufficiently large $n$, we have by (\ref{lem0100})
\begin{eqnarray*}
\epsilon b^{-m_n} &\leq& |b^{-m_n}u-x_0|\leq |b^{-m_n}u-y_n|+|y_n-x_0|\\
&\leq& b^{-m_n}(c+c_1)+b^{-m_n}\lambda_n\leq b^{-m_n}(c+c_1+\lambda+1).
\end{eqnarray*}
 This implies that
\begin{eqnarray*}
-\log|b^{-m_n}u-x_0|\sim m_n\log b.\label{lem0103}
\end{eqnarray*}
For sufficiently large $n$, it follows that
\begin{eqnarray*}
\int_{[x_n+c_1, x_n+c_1+c]\backslash E_n}h(\log u)du &=& \int_{[x_n+c_1, x_n+c_1+c]\backslash E_n}h(\log b^{-m_n}u)du\\
&=& \int_{[x_n+c_1, x_n+c_1+c]\backslash E_n}\frac{-1}{\log|b^{-m_n}u-x_0|}du\\
&\sim& \int_{[x_n+c_1, x_n+c_1+c]\backslash E_n}\frac{1}{m_n\log b}du\\
\end{eqnarray*}
As we have
\begin{eqnarray*}
&& c\geq \int_{[x_n+c_1, x_n+c_1+c]\backslash E_n}du\geq \int_{[x_n+c_1, x_n+c_1+c]}du-\int_{ E_n}du\geq c-2\epsilon,
\end{eqnarray*}
it follows that
\begin{eqnarray*}
&& (1-\epsilon)\cdot \frac{c-2\epsilon}{m_n\log b}\leq \int_{[x_n+c_1, x_n+c_1+c]\backslash E_n}h(\log u)du\leq (1+\epsilon)\cdot \frac{c}{m_n\log b}
\end{eqnarray*}
for sufficiently large $n$. Furthermore, we see from (\ref{lem0102}) that
\begin{eqnarray*}
\int_{[x_n+c_1, x_n+c_1+c]\cap E_n}h(\log u)du &=& \int_{[x_n+c_1, x_n+c_1+c]\cap E_n}\frac{-1}{\log|b^{-m_n}u-x_0|}du\\
&\leq& \frac{2}{m_n\log b} \int_{E_n}du\leq \frac{4\epsilon}{m_n\log b}.
\end{eqnarray*}
Hence, we obtain that
\begin{eqnarray*}
\int_{x_n+c_1}^{x_n+c_1+c}\phi(u)du &\sim& x_n^{-1-\alpha}\int_{x_n+c_1}^{x_n+c_1+c}h(\log u)du\\
&\sim& x_n^{-1-\alpha}\frac{c}{m_n\log b}, 
\end{eqnarray*}
so that (\ref{lem0101}) holds.\par
Case 2-2.\quad Suppose that $\lambda=\infty$. For $u$ with $x_n+c_1+\leq u\leq x_n+c_1+c$, we see from (\ref{lem0100}) that
\begin{eqnarray*}
|y_n-x_0|-(c+c_1)b^{-m_n}\leq |b^{-m_n}u-x_0|\leq |y_n-x_0|+(c+c_1)b^{-m_n},
\end{eqnarray*}
that is,
\begin{eqnarray*}
(1-(c+c_1)\lambda_n^{-1})|y_n-x_0|\leq |b^{-m_n}u-x_0|\leq (1+(c+c_1)\lambda_n^{-1})|y_n-x_0|.
\end{eqnarray*}
This implies that
\begin{eqnarray*}
\int_{x_n+c_1}^{x_n+c_1+c}\phi(u)du &\sim& x_n^{-1-\alpha}\int_{x_n+c_1}^{x_n+c_1+c}\frac{-1}{\log|b^{-m_n}u-x_0|}du\\
&\sim& x_n^{-1-\alpha}\frac{-c}{\log|y_n-x_0|},
\end{eqnarray*}
so we get (\ref{lem0101}). The lemma has been proved.  \qed

\vspace{2ex}

\begin{lem}
We have
\begin{eqnarray*}
\phi\otimes\phi(x)\sim 2M\int_{x}^{x+1}\phi(u)du=2M^2\mu((x,x+1]).
\end{eqnarray*}
\end{lem}

{\it Proof }\quad  Let $\{y_n\}$ be a sequence such that $1\leq y_n\leq b$ and $\lim_{n\to\infty}y_n=y$ for some $y\in[1,b]$. We put $x_n=b^{m_n}y_n$, where $m_n$ is a positive integer and $\lim_{n\to\infty}x_n=\infty$.  Now, we have
\begin{eqnarray*}
\phi\otimes\phi(x_n) &=& \int_1^{x_n-1}\phi(x_n-u)\phi(u)du\\
&=& 2\int_1^{2^{-1}x_n}\phi(x_n-u)\phi(u)du\\
&=& 2\left(\int_1^{(\log x_n)^\beta}+\int_{(\log x_n)^\beta}^{2^{-1}x_n}\right)\phi(x_n-u)\phi(u)du=:2(J_1+J_2).
\end{eqnarray*}
Here, we took $\beta$ satisfying $\alpha\beta>1$. Put $K:=\sup\{h(\log x) : 1\leq x\leq b \}$. Then, we have
\begin{eqnarray*}
J_2&\leq& K^2\int_{(\log x_n)^\beta}^{2^{-1}x_n}\frac{du}{u^{1+\alpha}(x_n-u)^{1+\alpha}}\leq K^2\left(\frac{2}{x_n}\right)^{1+\alpha}\cdot\alpha^{-1}(\log x_n)^{-\alpha\beta}.
\end{eqnarray*}
We consider the two cases where $y\not=x_0$ and $y=x_0$. \par
Case 1. \quad Suppose that $y\not=x_0$. If $1\leq u\leq (\log x_n)^\beta$, then
\begin{eqnarray*}
h(\log(x_n-u))=h(\log(y_n-b^{-m_n}u))\sim h(\log y).
\end{eqnarray*}
Hence, we obtain that 
\begin{eqnarray*}
J_1 &=& \int_1^{(\log x_n)^\beta}(x_n-u)^{-1-\alpha}u^{-1-\alpha}h(\log(x_n-u))h(\log u)du\\
&\sim& x_n^{-1-\alpha}\int_1^{(\log x_n)^\beta}u^{-1-\alpha}h(\log(x_n-u))h(\log u)du\\
&\sim& Mx_n^{-1-\alpha}h(\log y),
\end{eqnarray*}
so that
\begin{eqnarray*}
\phi\otimes\phi(x_n)=2(J_1+J_2)\sim 2J_1\sim 2 Mx_n^{-1-\alpha}h(\log y).
\end{eqnarray*}

Case 2.\quad Suppose that $y=x_0$. Put $\gamma_n:=b^{m_n}|y_n-x_0|(\log x_n)^{-\beta}$ and
\begin{eqnarray*}
E_n':=\{u : |y_n-x_0-b^{-m_n}u|\leq \epsilon b^{-m_n}\},
\end{eqnarray*}
where $0< \epsilon< 1$. It suffices that we consider the case where there exists a limit of $\gamma_n$, so we may put $\gamma:=\lim_{n\to\infty}\gamma_n$. This limit permits infinity. Furthermore, we divide $\gamma$ in the two cases where $\gamma<\infty$ and $\gamma=\infty$. 

Case 2-1.\quad Suppose that $0\leq \gamma<\infty$. Take sufficiently large $n$. Set
\begin{eqnarray*}
&& J_{11}':= \int_{[1,(\log x_n)^\beta]\backslash E_n'}u^{-1-\alpha}h(\log(x_n-u))h(\log u)du,\\
&& J_{12}':=\int_{[1,(\log x_n)^\beta]\cap E_n'}u^{-1-\alpha}h(\log(x_n-u))h(\log u)du.
\end{eqnarray*} 
 Let $u\in[1, (\log x_n)^\beta]\backslash E_n'$. We have
\begin{eqnarray*}
\epsilon b^{-m_n}\leq |y_n-x_0-b^{-m_n}u| &\leq& |y_n-x_0|+b^{-m_n}u\nonumber   \\ 
&\leq& (\gamma+2)b^{-m_n}(\log x_n)^\beta. \label{lem0201}
\end{eqnarray*}
This implies that
\begin{eqnarray*}
-\log|y_n-x_0-b^{-m_n}u|\sim m_n\log b.
\end{eqnarray*}
It follows that
\begin{eqnarray*}
J_{11}'&=& \int_{[1,(\log x_n)^\beta]\backslash E_n'}u^{-1-\alpha}h(\log(y_n-b^{-m_n}u))h(\log u)du\\
&=& \int_{[1,(\log x_n)^\beta]\backslash E_n'}u^{-1-\alpha}h(\log u)\frac{-1}{\log|y_n-x_0-b^{-m_n}u|}du\\
&\sim& \frac{1}{m_n\log b} \int_{[1,(\log x_n)^\beta]\backslash E_n'}u^{-1-\alpha}h(\log u)du.
\end{eqnarray*}
Here, we see that, for sufficiently large $n$, 
\begin{eqnarray*}
M-\epsilon-2\epsilon K\leq \int_{[1,(\log x_n)^\beta]\backslash E_n'}u^{-1-\alpha}h(\log u)du\leq M,
\end{eqnarray*}
and thereby
\begin{eqnarray*}
(1-\epsilon)\frac{M-\epsilon-2\epsilon K}{m_n\log b}\leq J_{11}'\leq (1+\epsilon)\frac{M}{m_n\log b}.
\end{eqnarray*}
Let $u\in E_n'$. Then, we have 
\begin{eqnarray*}
h(\log(x_n-u)) &=& h(\log(y_n-b^{-m_n}u))\\
&=& \frac{-1}{\log|y_n-x_0-b^{-m_n}u|}\leq \frac{2}{m_n\log b}.
\end{eqnarray*}
Hence, we see that
\begin{eqnarray*}
J'_{12} &\leq& \frac{2}{m_n\log b}\int_{[1,(\log x_n)^\beta]\cap E_n'}u^{-\alpha -1}h(\log u)du   \leq  \frac{4K\epsilon}{m_n\log b}.
\end{eqnarray*}
We consequently obtain that
\begin{eqnarray*}
J_1\sim x_n^{-1-\alpha}(J_{11}'+J_{12}')\sim \frac{Mx_n^{-1-\alpha}}{m_n\log b},
\end{eqnarray*}
so that
\begin{eqnarray*}
\phi\otimes\phi(x_n)=2(J_1+J_2)\sim 2J_1\sim \frac{2Mx_n^{-1-\alpha}}{m_n\log b}.
\end{eqnarray*}

Case 2-2.\quad Suppose that $\gamma=\infty$. Note that $[1, (\log x_n)^\beta]\cap E_n'$ is empty for sufficiently large $n$. Let $1\leq u\leq (\log x_n)^\beta$. Since
\begin{eqnarray*}
&&|y_n-x_0|(1-\gamma_n^{-1})\leq |y_n-x_0-b^{-m_n}u|\leq |y_n-x_0|(1+\gamma_n^{-1}),
\end{eqnarray*}
 we see that
\begin{eqnarray*}
\log |y_n-x_0-b^{-m_n}u|\sim \log |y_n-x_0|.
\end{eqnarray*}
This yields that
\begin{eqnarray*}
J_1 &\sim& x_n^{-1-\alpha}\int_{[1,(\log x_n)^\beta]}u^{-1-\alpha}h(\log u)\cdot\frac{-1}{\log|y_n-x_0-b^{-m_n}u|}du\\
 &\sim& \frac{-M}{\log|y_n-x_0|}x^{-1-\alpha}_n.
\end{eqnarray*}
For sufficiently large $n$, we have
\begin{eqnarray*}
J_2\times x_n^{1+\alpha}(-\log|y_n-x_0|)&\leq& \frac{2^{1+\alpha}K^2}{\alpha}\cdot \frac{-\log|y_n-x_0|}{(\log x_n)^{\alpha\beta}}\\
&=& \frac{2^{1+\alpha}K^2}{\alpha}\cdot \frac{-\log\gamma_n+m_n\log b-\log(\log x_n)^\beta}{(\log x_n)^{\alpha\beta}}\\
&\leq& \frac{2^{1+\alpha}K^2}{\alpha}\cdot \frac{m_n\log b}{(\log x_n)^{\alpha\beta}},
\end{eqnarray*}
so that ${\displaystyle \lim_{n\to\infty}J_2/J_1=0}$. We consequently obtain that
\begin{eqnarray*}
&& \phi\otimes\phi(x_n)=2(J_1+J_2)\sim 2J_1\sim 2x_n^{-1-\alpha}\frac{-M}{\log|y_n-x_0|}.
\end{eqnarray*}

Combining the above calculations with the proof of Lemma 3.1, we reach the following: If $y\not=x_0$, then
\begin{eqnarray*}
&& \phi\otimes\phi(x_n)\sim 2Mx_n^{-1-\alpha}h(\log y)\sim 2M\int_{x_n}^{x_n+1}\phi(u)du.
\end{eqnarray*}
Suppose that $y=x_0$. Recall $\lambda$ in the proof of Lemma 3.1. 
If $0\leq \gamma<\infty$ and $\lambda=\infty$, then we have $-\log|y_n-x_0|\sim m_n\log b$.
Hence,
\begin{eqnarray*}
\phi\otimes\phi(x_n)&\sim& 2M\frac{x_n^{-1-\alpha}}{m_n\log b}\\
&\sim& 2M\frac{-x_n^{-1-\alpha}}{\log|y_n-x_0|}\sim 2M\int_{x_n}^{x_n+1}\phi(u)du.
\end{eqnarray*}
If $0\leq \gamma<\infty$ and $0\leq \lambda<\infty$, then
\begin{eqnarray*}
\phi\otimes\phi(x_n)&\sim& 2M\frac{x_n^{-1-\alpha}}{m_n\log b}\sim 2M\int_{x_n}^{x_n+1}\phi(u)du.
\end{eqnarray*}
If $\gamma=\infty$, then $\lambda=\infty$ and
\begin{eqnarray*}
&& \phi\otimes\phi(x_n)\sim 2M\frac{-x_n^{-1-\alpha}}{\log|y_n-x_0|}\sim 2M\int_{x_n}^{x_n+1}\phi(u)du.
\end{eqnarray*}
The lemma has been proved. \qed

\vspace{2ex}

{\it Proof of Theorem 1.1 }\quad We have $\mu\in\mathcal{L}_{loc}$ by Lemma 3.1. It follows from Lemma 3.2 that
\begin{eqnarray*}
\mu*\mu((x, x+1])&=&M^{-2}\int_x^{x+1}\phi\otimes\phi(u)du\\ \nonumber
&\sim& 2\int_x^{x+1}\mu((u,u+1])du\sim 2\mu((x,x+1]).
\end{eqnarray*}
Let $c>0$. Furthermore, we see from $\mu\in\mathcal{L}_{loc}$ and (ii) of Lemma 2.2 that
\begin{eqnarray*}
\mu*\mu((x, x+c])\sim c\mu*\mu((x,x+1])\quad \mbox{and}\quad \mu((x, x+c])\sim c\mu((x,x+1]).\nonumber
\end{eqnarray*}
Hence, we get
\begin{eqnarray*}
\mu*\mu((x, x+c])\sim 2\mu((x,x+c]),
\end{eqnarray*}
and thereby $\mu\in\mathcal{S}_{loc}$. Thus, $\mu((x-1,x])\in \mathcal{S}_{d}$ by (i) of Proposition 1.1. Since we see  that 
\begin{eqnarray*}
\phi\otimes\phi(x)\sim 2M\int_{x}^{x+1}\phi(u)du=2M^2\mu((x,x+1]),
\end{eqnarray*}
we have $\mu^{2*} \in \mathcal{S}_{ac}$ by (ii) of Lemma 2.1.
However, we have $\mu\not\in\mathcal{UL}_{loc}$ because, for $c= b^{-m(n)}$ with $m(n) \in \N$, we see that as $n \to \infty$
\begin{equation}
c^{-1}\int_{b^nx_0}^{b^nx_0 +c}M^{-1}\phi(u)du \sim \frac{M^{-1}b^{-(\alpha+1)n}x_0^{-\alpha-1}}{(m(n)+n)\log b}.\nonumber
\end{equation}
The above relation implies that the convergence of the definition of 
the class $\mathcal{UL}_{loc}$ fails to satisfy uniformity.  Since $\mathcal{US}_{loc}\subset \mathcal{UL}_{loc}$, the theorem has been proved.  \qed

\medskip

{\it Proof of Corollary 1.1 } \q   Proofs of assertions (i) and (ii) are clear from Theorem 1.1. We find from the proof of Theorem 1.1 that $\mu \notin \mathcal{UL}_{loc}$ but $\mu^{2*} \in \mathcal{S}_{ac}
$. Since $\mathcal{S}_{ac}\subset \mathcal{L}_{ac}\subset \mathcal{UL}_{loc}$, assertions (iii) and (iv) are true. \qed

\medskip

Choose $x_1$ and $x_2$ satisfying that $1 < x_0 < x_0 +x_1 <x_0 +x_2 <b$. Let $\{n_k\}_{k=1}^{\infty}$ be an  increasing sequence of positive integers satisfying $\sum_{k=1}^{\infty}1/\sqrt{n_k}=1.$ Let $B_k:=(-b^{n_k}x_2,-b^{n_k}x_1]$ and $D_k:=(b^{n_k}x_0,b^{n_k}x_0+1]$ for $k \in \N$. Choose a distribution $\mu_1$ satisfying that
$\mu_1(B_k)= 1/\sqrt{n_k}$  for all $k \in \N$ and $\mu_1((\cup_{k=1}^{\infty}B_k)^c)=0.$

\begin{lem} We have, for $ c \in \R$,
\begin{equation}
\lim_{k \to \infty}\frac{\mu*\mu_1(D_k+c)}{\mu(D_k)} = \infty.\nonumber
\end{equation}
\end{lem}
{\it Proof } \q  We have, uniformly in $v \in [x_1,x_2]$, 
\begin{equation}
\mu((b^n(x_0 +v),b^n(x_0+v)+1])\sim  M^{-1} b^{-(\alpha +1)n}(x_0 +v)^{-\alpha -1}h(\log(x_0+v))\nonumber
\end{equation}
and
\begin{equation}
\mu((b^nx_0,b^nx_0+1]) \sim  M^{-1}\frac{b^{-(\alpha +1)n}x_0^{-\alpha -1}}{n \log b}.\nonumber
\end{equation}
Thus,  there exists $c_1 >0$ such that $c_1$ does not depend on $ v \in [x_1, x_2]$ and that 
\begin{equation}
\liminf_{n \to \infty}\frac{\mu((b^n(x_0 +v),b^n(x_0+v)+1])}{n\mu((b^nx_0,b^nx_0+1])} \ge c_1.\nonumber
\end{equation}
Hence, we obtain from Lemma 3.1  that
\begin{equation}
\begin{split}
\liminf_{k \to \infty}\frac{\mu*\mu_1(D_k+c)}{\mu(D_k)} &\ge  \liminf_{k \to \infty}\int_{B_k}\frac{\mu(D_k-u+c)}{\mu(D_k)}\mu_1(du)\\ \nonumber
& =  \liminf_{k \to \infty}\int_{B_k}\frac{\mu(D_k-u)}{\mu(D_k)}\mu_1(du)\\ \nonumber
&\ge c_1 \liminf_{k \to \infty}\frac{ n_k}{\sqrt{n_k}}= \infty.\nonumber
\end{split}
\end{equation}
Thus, we have proved the lemma. \qed

\medskip

{\it Proof of Theorem 1.2 } \q  Define  distributions $\rho_1 $ and $\rho_2 $ as 
\begin{equation}
\rho_1(dx):= 2^{-1}\delta_{0}(dx) + 2^{-1}\mu(dx), \q \rho_2(dx):= 2^{-1}\mu_1(dx) + 2^{-1}\mu(dx).\nonumber
\end{equation}
Thus, $\rho_1 \in \mathcal{S}_{loc}$ by Theorem 1.1 and (iii) of Lemma 2.2. Let $\rho(dx):= f(x)dx$, where $f(x)$ is continuous with compact support in $[0,1]$. Define  distributions $p_1(x)dx$ and $p_2(x)dx$ as 
\begin{equation}
p_1(x)dx:= \rho*\rho_1(dx)=2^{-1}f(x)dx + 2^{-1}\rho*\mu(dx) \nonumber
\end{equation}
and
\begin{equation}
  p_2(x)dx:= \rho*\rho_2(dx)=2^{-1}\rho*\mu_1(dx) + 2^{-1}\rho*\mu(dx). \nonumber
\end{equation}
Then, we find that $p_1(x) =p_2(x)$ for all sufficiently large $x>0$ and $p_1(x)\in \mathcal{S}_{d}$ by (ii) of Proposition 1.1. 
 We establish from  Lemma 3.3 and Fatou's lemma that
\begin{equation}
\begin{split}
&\liminf_{k \to \infty}\frac{\int_{D_k}p_2\otimes p_2(x)dx}{\int_{D_k}p_2(x)dx}\\\nonumber
& \ge \liminf_{k \to \infty}\frac{\int_0^2\mu*\mu_1(D_k-u)f^{2\otimes}(u)du}{\int_0^1\mu(D_k-u)f(u)du} \\\nonumber
&\ge  \int_0^2\liminf_{k \to \infty}\frac{\mu*\mu_1(D_k-u)}{\mu(D_k)} f^{2\otimes}(u)du= \infty.\nonumber
\end{split}
\end{equation}
Thus, we conclude that $p_2(x)\notin \mathcal{S}_{d}$.
\qed


\section{A remark on the closure under convolution roots }

The  tail of a measure $\xi$ on $\R$ is denoted by  $\bar \xi(x)$, that is, $\bar \xi(x): = \xi((x,\infty))$ for $x \in \R$. Let $\gamma \in \R$. The $\gamma$-exponential moment of $\xi$ is denoted by $\widehat \xi(\gamma)$, namely, $\widehat \xi(\gamma):= \int_{-\infty}^{\infty}e^{\gamma x}\xi(dx)$. 

\begin{defn}\quad Let $\gamma \ge 0$.

 (i) A distribution $\rho$ on  $\R$ is said to belong to the class $\mcal{ L}(\gamma)$ if $\bar \rho(x) >0$ for every $x \in  \R$ and if
\begin{equation}\label{2.1}
 \bar \rho(x+a)\sim  e^{-\gamma a}\bar \rho(x) \q \mbox{ for every  }\q  a \in \R. \nonumber
\end{equation}

(ii)  A distribution $\rho$ on $\mathbb{R}$ belongs to the class ${\mcal S}(\gamma)$ if $\rho \in {\mcal L}(\gamma)$ with $\widehat \rho(\gamma)< \infty$ and if
\begin{equation}\label{2.2}
 \overline{\rho*\rho}(x) \sim  2\widehat \rho(\gamma)\bar \rho(x). \nonumber
\end{equation}

(iii) Let $\gamma_1 \in \R$. A distribution $\rho$ on $\mathbb{R}$ belongs to the class ${\mcal M}(\gamma_1)$ if $\widehat \rho(\gamma_1)< \infty$. 

\end{defn}
\medskip

The convolution closure problem on the class $\mathcal{S}(\gamma)$  with  $\gamma \ge 0$ is negatively solved by Leslie \cite{le89} for $\gamma = 0$ and by 
Kl\"uppelberg and  Villasenor \cite{kv91} for  $\gamma > 0$. The same problem on the class $\mathcal{S}_{d}$ is also negatively solved by Kl\"uppelberg and  Villasenor \cite{kv91}.  On the other hand,  the fact that the class $\mathcal{S}(0)$ of subexponential distributions is closed under convolution roots is proved by Embrechts et al.\ \cite{egv79} in the one-sided case and  by Watanabe \cite{wa08} in the two-sided case. Embrechts and Goldie conjecture that ${\mcal L}(\gamma)$  with  $\gamma \ge 0$  and $\mathcal{S}(\gamma)$  with  $\gamma > 0$  are  closed under convolution roots in \cite{eg80,eg82}, respectively. They also prove  in \cite{eg82} that if  $\mathcal{L}(\gamma)\cap{\mcal P}_+$ with $\gamma > 0$ is  closed under convolution roots, then $\mathcal{S}(\gamma)\cap{\mcal P}_+$  with $\gamma > 0$ is  closed under convolution roots. However, Shimura and Watanabe \cite{shw} prove that the class ${\mcal L}(\gamma)$  with  $\gamma \ge 0$ is not closed under convolution roots, and we find that Xu et al. \cite{xu} show the same conclusion in the case $\gamma=0$. Pakes \cite{pa04} and Watanabe \cite{wa08} show that $\mathcal{S}(\gamma)$ with $\gamma > 0$ is closed under convolution roots in the class of infinitely divisible distributions on $\R$. It is still open whether the class $\mathcal{S}(\gamma)$  with  $\gamma >0$ is closed under convolution roots. Shimura and Watanabe \cite{shw05} show that the class $\mathcal{OS}$ is not closed under convolution roots. Watanabe and Yamamuro  \cite{wy10a} pointed out that $\mathcal{OS}$ is closed under convolution roots in the class of infinitely divisible distributions.  

\vspace{1ex}

\medskip

Let $\gamma \in \R$. For $ \mu \in  {\mcal M}(\gamma)$, we define the {\it exponential tilt} $\mu_{\langle\gamma\rangle}$ of $\mu$ as
\begin{equation}
\mu_{\langle\gamma\rangle}(dx):= \frac{1}{\widehat\mu(\gamma)}e^{\gamma x}\mu(dx). \nonumber
\end{equation}
Exponential tilts preserve convolutions, that is, $(\mu*\rho)_{\langle\gamma\rangle}=\mu_{\langle\gamma\rangle}*\rho_{\langle\gamma\rangle}$ for distributions $\mu, \rho \in {\mcal M}(\gamma)$. Let  ${\mcal C}$ be a distribution class. For a  class ${\mcal C}\subset {\mcal M}(\gamma)$, we define the class ${\frak E}_{\gamma}( {\mcal C})$ by
\begin{equation}
{\frak E}_{\gamma}( {\mcal C}):=\{\mu_{\langle\gamma\rangle} : \mu \in {\mcal C}\}. \nonumber
\end{equation}
It is obvious that ${\frak E}_{\gamma}(  {\mcal M}(\gamma))=  {\mcal M}(-\gamma) $ and that $(\mu_{\langle\gamma\rangle})_{\langle-\gamma\rangle}= \mu$ for $\mu \in {\mcal M}(\gamma)$.   The class ${\frak E}_{\gamma}( {\mcal S}(\gamma))$ is determined  by Watanabe and Yamamuro as follows. 
Analogous result is found in Theorem 2.1  of Kl\"uppelberg \cite{kl89}.

\begin{lem} (Theorem 2.1 of \cite{wy10})
Let $\gamma >0$.

(i) We have ${\frak E}_{\gamma}( {\mcal L}(\gamma)\cap {\mcal M}(\gamma))= {\mcal L}_{loc}\cap {\mcal M}(-\gamma) $ and hence ${\frak E}_{\gamma}( {\mcal L}(\gamma)\cap {\mcal M}(\gamma)\cap{\mcal P}_+)= {\mcal L}_{loc}\cap{\mcal P}_+ $. Moreover, if $ \rho \in {\mcal L}(\gamma)\cap {\mcal M}(\gamma)$, then we have 
\begin{equation}
\rho_{\langle\gamma\rangle}((x,x+c]) \sim \frac{c\gamma}{\widehat\rho(\gamma)}e^{\gamma x}\bar\rho(x) \mbox{ for all } c>0. \nonumber
\end{equation}

(ii)  We have ${\frak E}_{\gamma}( {\mcal S}(\gamma))= {\mcal S}_{loc}\cap {\mcal M}(-\gamma) $ and thereby  ${\frak E}_{\gamma}( {\mcal S}(\gamma)\cap{\mcal P}_+)= {\mcal S}_{loc}\cap{\mcal P}_+ $.
\end{lem} 
 
\medskip

Finally, we present a remark on the closure under convolution roots for the three classes ${\mcal S}(\gamma)\cap{\mcal P}_+$, ${\mcal S}_{loc}\cap {\mcal P}_+$, and ${\mcal S}_{ac}\cap {\mcal P}_+$.

\begin{prop} The following are equivalent:

(1)  The class $\mathcal{S}(\gamma)\cap \mathcal{P}_{+}$ with $\gamma >0$ is closed under convolution roots. 

(2) The class $\mathcal{S}_{loc}\cap \mathcal{P}_{+}$ is closed under convolution roots.

(3) Let $\mu$ be a distribution on $\R_+$ and let $p_c(x):=c^{-1}\mu((x-c,x]) $ for  $c>0$. Then,   $\{p_c^{n\otimes}(x) : c>0\} \subset \mathcal{S}_{d}$   for some $n \in \N$  implies $\{p_c(x) : c>0\} \subset  \mathcal{S}_{d}$.  
\end{prop}

{\it Proof } \q   Proof of the equivalence between (1) and (2)  is due to Lemma 4.1. Let $n \ge 2$. Suppose that (2) holds and,  for some $n$,   $p_c^{n\otimes}(x)\in \mathcal{S}_{d}$ for every $c>0$. Let $f_c(x) =c^{-1}1_{[0,c)}(x)$. We have $p_c^{n\otimes}(x)dx= ((f_c(x)dx)*\mu)^{n*}\in \mathcal{S}_{loc}$. We see from assertion (2) that $(f_c(x)dx)*\mu\in \mathcal{S}_{loc}$ and hence, by (iii) of Proposition 1.1, we have $\mu \in \mathcal{S}_{loc}$, that is, $p_c(x) \in \mathcal{S}_{d}$ for every $c>0$ by (i) of Proposition 1.1.
Conversely, suppose that (3) holds and $\mu^{n*}\in \mathcal{S}_{loc}$. Note that $f_c^{n\otimes}(x)$ is continuous with compact support in $\R_+$. Thus, we see from (ii) of Proposition 1.1 that $p_c^{n\otimes}(x)= \int_{0-}^{x+}f_c^{n\otimes}(x-u)\mu^{n*}(du)\in\mathcal{S}_{d}$ for every $c>0$. We obtain from assertion (3) that  $p_c(x)\in \mathcal{S}_{d}$ for every $c>0$, that is, $\mu \in \mathcal{S}_{loc}$ by (i) of Proposition 1.1. \qed


\end{document}